\documentclass{amsart}

\usepackage{amsmath, amsthm, amsfonts, amssymb}
\usepackage{mathtools}
\usepackage{graphicx}
\usepackage{MnSymbol} 
\usepackage{longtable}
\usepackage{caption}

\newtheorem{thm}{Theorem}[section]
\newtheorem{cor}[thm]{Corollary}
\newtheorem{lem}[thm]{Lemma}
\newtheorem{prop}[thm]{Proposition}
\theoremstyle{definition}
\newtheorem{defn}[thm]{Definition}
\newtheorem*{strat}{Strategy}
\theoremstyle{remark}
\newtheorem{rem}[thm]{Remark}

\def\ZZ{\mathbb{Z}}
\def\NN{\mathbb{N}}
\def\FF{\mathbb{F}}
\def\PP{\mathbb{P}}
\def\|{\left | \right.}

\def\CG{{\mathcal G}}

\def\CG{{\mathcal G}}

\def\CL{{\mathcal L}}

\def\CP{{\mathcal P}}

\def\CT{{\mathcal T}}

\makeatletter
\newcommand{\bigperp}{%
  \mathop{\mathpalette\bigp@rp\relax}%
  \displaylimits
}
\newcommand{\bigp@rp}[2]{%
  \vcenter{
    \m@th\hbox{\scalebox{\ifx#1\displaystyle2.1\else1.5\fi}{$#1\perp$}}
  }%
}
\makeatother

\newcommand{\sA}{\mathsf A}  

\newcommand{\sD}{\mathsf D}

\newcommand{\sR}{\mathsf R}

\newcommand{\xsA}[1]{\hspace{0.5mm}^{#1}\hspace*{-0.5mm}\mathsf A}
\newcommand{\xsB}[1]{\hspace{0.5mm}^{#1}\hspace*{-0.2mm}\mathsf B}
\newcommand{\xsC}[1]{\hspace{0.5mm}^{#1}\hspace*{-0.2mm}\mathsf C}
\newcommand{\xsD}[1]{\hspace{0.5mm}^{#1}\hspace*{-0.2mm}\mathsf D}
\newcommand{\xsE}[1]{\hspace{0.5mm}^{#1}\hspace*{-0.2mm}\mathsf E}
\newcommand{\xsF}[1]{\hspace{0.5mm}^{#1}\hspace*{-0.2mm}\mathsf F}
\newcommand{\xsG}[1]{\hspace{0.5mm}^{#1}\hspace*{-0.2mm}\mathsf G}

\newcommand{\xsR}[1]{\hspace{0.5mm}^{#1}\hspace*{-0.2mm}\mathsf R}

\DeclareMathOperator{\Quot}{Quot}
\DeclareMathOperator{\std}{std}
\DeclareMathOperator{\m}{m}
\DeclareMathOperator{\cp}{m_{p}}
\DeclareMathOperator{\mref}{m_{ref}}
\DeclareMathOperator{\M}{M}
\DeclareMathOperator{\N}{N}
\DeclareMathOperator{\Mref}{M_{ref}}
\DeclareMathOperator{\Nref}{N_{ref}}
\DeclareMathOperator{\D}{D}
\DeclareMathOperator{\E}{E}

\begin{document}

\title{Totally-Reflective Genera of Integral Lattices}
\author{Ivica Turkalj}

\begin{abstract} 
 In this paper we give a complete classification of totally-reflective, primitive genera in dimension 3 and 4. Our method breaks up into two parts. 
 The first part consists of classifying the square free, totally-reflective, primitive genera by calculating strong bounds on the prime factors of the determinant of positive definite
 quadratic forms (lattices) with this property. 
 We achieve these bounds by combining the Minkowski-Siegel mass formula with the combinatorial classification of reflective lattices accomplished by Scharlau \& Blaschke. In a second part,
 we use a lattice transformation that goes back to Watson, to generate all totally-reflective, primitive genera when starting with the square free case. 
\end{abstract}

\maketitle

\section{Introduction}
The investigation of totally-reflective genera came up with the problem of classifying the cofinite, arithmetic reflection groups on the hyperbolic space. 
These groups occur as Weyl groups of Lorentzian lattices with a fundamental polyhedron of finite volume. 
Vinberg showed in \cite{vin72} that a necessary condition for (certain) Lorentzian lattices to induce a cofinite, arithmetic reflection group,
is that a naturally associated, positive definite genus is totally-reflective. In the subsequent work \cite{vin81}, he proved 
that totally-reflective genera only appear in dimension $< 30$. It was unknown whether this bound is sharp since the largest known dimension of such genera was $20$.
An example was found by Borcherds in \cite{bor87}, \S$8$, Example $5$. 
By a more detailed investigation of the existence of non reflective lattices in high dimensional genera, Esselmann proved in \cite{ess96} 
that $20$ is actually the largest dimension in which a totally-reflective genus can exists. Furthermore, it is well known that there are only finitely many totally-reflective, primitive
genera in any fixed
dimension (cf. \cite{sw92}, Theorem 1.4). An explicit classification has been carried out only in dimenson 3 for the special case of square free genera, cf. \cite{wal93} (which we reproduce here). 
In this work, all totally-reflective, primitive genera in dimension $3$ and $4$ have been found. Also, the related problem of classifying the hyperbolic reflection groups is outlined. 
We give a brief overview of our methods.
\begin{strat} Our goal of classifying all totally-reflective genera in dimension $3$ and $4$ is achieved as follows:
 \begin{itemize}
  \item[Step 1:] Let $L$ be a strongly square free, totally-reflective lattice with $\dim L = 4$ (resp. $\dim = 3$). Hence $\det L$ is of the form $\det L = p_1^2 \cdots p_r^2 \cdot q_1 \cdots q_s$
                 (resp. $r=0$). 
                 Using the mass formula, we prove that $r \leq 9$ and $s \leq 8-r$ (resp. $s \leq 10)$ (section 3.1).
  \item[Step 2:] By applying the combinatorial classification of Scharlau \& Blaschke, we prove that there are bounds
                 $\overline{p_i}$ and $\overline{q_j}$ (one for every prime factor) depending only 
                 on the number of prime factors,
                 such that $p_i \leq \overline{p_i}$ and $q_j \leq \overline{q_j}$. Thus the number of local invariants that need to be taken into account
                 is effectively bounded and the enumeration is computationally feasible (section 3.2)
  \item[Step 3:] After finishing the strongly square free classification, we obtain all square free, totally-reflective, primitive genera by partial dualization (section 4.1).
  \item[Step 4:] The last step consists in dropping the assumption ``square free'' by determining the pre-images of square free genera under the Watson transformation (section 4.2).

 \end{itemize}
\end{strat}

\section{Background}
\subsection{Integral lattices}
Let $R$ be a prinicipal ideal domain and $K:=\Quot(R)$ its quotient field.
A \emph{lattice} over $R$ is a pair $(L,b)$, where $L$ is a free $R$-module of finite rank 
and $b:L \times L \longrightarrow K$ a symmetric bilinear form. As usual, $V:=L \otimes_R K$ means the enveloping $K$-space of $L$. By $O(L)$ we denote the isometry group of $(L,b)$.
The \emph{determinant} of $L$ is the determinant of any Gram matrix of $(L,b)$, which is well defined modulo squares of units in $R$.

We say that $L$ is \emph{integral}, if $b(L,L) \subseteq R$.
An $R$-lattice is called \emph{even} if $b(x,x) \in 2R$ for all $x \in L$, and \emph{odd} otherwise.
By ${}^{\alpha}\!{L}$ we mean the lattice $(L,\alpha b)$ obtained by scaling the bilinear form, where $\alpha \in K$.
An integral lattice is said to be \emph{primitive} if it is not the scaled version of another intergral lattice.
We denote by $L^{\#}$ the \emph{dual lattice} of $L$ which is defined as
\begin{align*}
 L^{\#} := \left\{ v \in V \| \forall x \in L : b(x,v) \in R \right\}.
\end{align*}
Clearly, $L$ is integral iff $L \subseteq L^{\#}$. For an integral lattice, the group $L^{\#}/L$ has order $\det L$ and is called the \emph{discriminant group} of $L$.
A lattice is \emph{unimodular} if $L = L^{\#}$.
More generally, for $\alpha \in K$, we say $L$ is \emph{$\alpha$-modular} if $L = \alpha L^{\#}$.
It is easy to see that any $\alpha$-modular lattice $K$ can be written as $K={}^{\alpha}\!{L}$, where $L$ is unimodular.

Let $\PP$ be the set of all rational primes. Two $\ZZ$-lattices $L_1, L_2$ are in the same \emph{genus} if they become isometric over all completions $\ZZ_p$:
\begin{align*}
  L_1 \otimes_{\ZZ} \ZZ_p \cong L_2 \otimes_{\ZZ} \ZZ_p, \ \text{for all $p  \in \PP \cup \{ \infty \}$}.
\end{align*}
It is well known that any genus consists of finitely many isometry classes. We write $\mathcal{G}(L)$ for the set of all isometry classes in the genus of $L$ and 
define $h(L):=\# \mathcal{G}(L)$ as the \emph{class number} of $L$.

Recall that every $\ZZ$-lattice $L$ locally posseses a \emph{Jordan decomposition},
\begin{align}
 L \otimes_{\ZZ} \ZZ_p = {}^{p^{-r}}\!{L_{-r}} \perp \dots \perp {}^{p^{-1}}\!{L_{-1}} \perp L_0 \perp {}^{p}\!{L_1} \perp \dots \perp {}^{p^r}\!{L_r},
\end{align}
where all $L_i$ are unimodular (possibly zero-dimensional). For $p \neq 2$ the Jordan decomposition is unique up to isometry. 
Unfortunately this is not true for $p = 2$. At least the following data remains invariant: $\dim L_i$, $2^i$ and the property even/odd for every $L_i$. 
We refer to $L$ as \emph{square free} if the Jordan decomposition of $L$ is of the form $L_0 \perp {}^{p}\!{L_1}$ at every prime $p$, and as \emph{strongly square free} 
if additionally $\dim L_0 \geq \dim L_1$ holds. 

We remind the reader of the \emph{genus symbol} as introduced in \cite{cs99}, Chapter $15$. 
This symbol is a list of local symbols for each prime $p$ dividing $2\det L$. Assuming a Jordan decomposition as above, the local symbol at the prime $p \neq 2$ is the formal product
\begin{align*}
 \text{$\prod_{i=-r}^{s} (p^{i})^{\varepsilon_i,n_i}$, where $\varepsilon_i := \left(\frac{\det L_i}{p}\right)$ and $n_i:=\dim L_i$.}
\end{align*}
Due to the lack of uniqueness, the $2$-adic symbol is more complicated and two other invariants have to be taken into account (oddity and parity).
We refer to \cite{cs99} for a more detailed investigation.

\subsection{The mass formula}
From now on all $\ZZ$-lattices $L$ are assumed to be positive definite (except the Lorentzian lattices in section $5$). Hence $O(L)$ is a finite group and the following definition makes sense.

\begin{defn}Let $L$ be a $\ZZ$-lattice. The \emph{mass} of $L$ is defined as
 \begin{align*}
  \m(L):= \sum_{M \in \mathcal{G}(L)}\frac{1}{|O(M)|}.
 \end{align*}
\end{defn}
An indispensable tool for our investigation is the Minkowski-Siegel mass formula
which relates the mass of a lattice to local quantities, which can be derived from the genus symbol. Originally the mass formula is stated in terms of $p$-adic densities, 
cf. \cite{sie35} for further details.
However, for computational reasons, we found the approach of Conway and Sloane more suitable, cf. \cite{cs88}. Below we outline the fundamental aspects. 
Let $n:=\dim L, s:=\lceil \frac{n}{2} \rceil$ and $D:=(-1)^s\det L$.
A helpful notion is the so-called \emph{standard mass},
\begin{align*}
 \std(n,D):= 2\pi^{-n(n+1)/4} \cdot \prod_{j=1}^{n}\Gamma(\tfrac{j}{2}) \cdot \zeta(2)\zeta(4) \cdots \zeta(2s-2)\zeta_D(s),
\end{align*}
where $\Gamma$ denotes the Gamma function, $\zeta$ the Riemann zeta function and $\zeta_D$ refers to the $L$-function
\begin{align*}
 \zeta_D(s)=
  \begin{cases}
   \prod_{p \in \PP}\left(1-\left(\tfrac{D}{p}\right)\frac{1}{p^s}\right)^{-1}, & n\ \text{even}, \\
   1,                                                                           & n\ \text{odd}.
  \end{cases}
\end{align*}
The actual mass of $L$ is gained from the standard mass by multiplying with certain correction factors, one for every prime $p$ dividing $2\det L$. 
Unlike the standard mass, these correction factors depend on the local structure of $L$.

\begin{prop}[cf. \cite{cs88}]
 Let $L$ be a $\ZZ$-lattice with Jordan decompositions as in $(1)$. Then
 \begin{align*}
  \m(L) = \std(n,D) \cdot \prod_{p \mid 2\det L} \left( \cp(L) \cdot 2 \prod_{j=2}^{s}\left(1-p^{2-2j}\right) \right),
 \end{align*}
 where for $p \neq 2$
 \begin{align*}
  \cp(L) =\prod_{\substack{i \in \ZZ,\\ \dim L_i \neq 0}}\left(\tfrac{1}{2}\left(1+\varepsilon p^{-s_i}\right)^{-1} \cdot \prod_{i=2}^{s_i}\left(1-p^{2-2i}\right)^{-1}\right)
           \cdot \prod_{\substack{k,l \in \ZZ,\\ k<l}} p^{\tfrac{1}{2}(l-k)n_ln_k}.
 \end{align*}
\end{prop}
In the above proposition, $\varepsilon = 0$ if $\dim L_i$ is odd, and $\varepsilon \in \left\{1,-1 \right\}$ if $\dim L_i$ is even. The exact value of $\varepsilon$
depends on the species of the orthogonal group $O_{n_i}(p)$ over $\FF_p$, which can be read off the genus symbol. Again, the case $p=2$ is more complicated, for which we refer to \cite{cs88}.

\subsection{Reflective lattices and totally-reflective genera}
For a nonzero vector $v \in L$, the reflection
\begin{align*}
 s_v: x \longmapsto x - \frac{2b(x,v)}{b(v,v)}v
\end{align*}
is an isometry of $V$. We call $v$ a \emph{root} of $L$ if $v$ is primitive (that is $v/m \notin L$ for all integers $m > 0$) and $s_v(L) = L$. 
The set $R(L)$ of all roots of $L$ is a root system in the usual sense of Lie theory.
The subgroup $W(L) \leqslant O(L)$ generated by all reflections $s_v$, with $v \in R(L)$, is called \emph{Weyl group} of $L$.
Besides the combinatorial structure induced by the Dynkin diagram, the root system of a lattice inherits the quadratic form,
thus $R(L)$ decomposes into scaled irreducible components $\xsA{\alpha}_n, \xsB{\alpha}_n, \xsC{\alpha}_n, \xsD{\alpha}_n, 
\xsE{\alpha}_6, \xsE{\alpha}_7, \xsE{\alpha}_8, \xsF{\alpha}_4, \xsG{\alpha}_2$. We mention \cite{bou08} as a general reference to the theory of root systems.

A positive definite lattice is called reflective if, roughly spoken, it is ``almost'' a root lattice (it is a root lattice up to finte index). More precisely:
\begin{defn}
 The positive definite, integral lattice $L$ is called \emph{reflective} if its root system $R(L)$ generates a sublattice of the same rank.
\end{defn}
Certainly, $L$ is reflective iff $W(L)$ has no nonzero fixed vectors (while acting on $V$). 

\begin{rem}
Crucial to our investigations is the work of Scharlau \& Blaschke. 
They classified all indecomposable, reflective lattices in low dimensions by pairs $(\sR, \CL)$, 
where $\sR$ is a scaled root system and $\CL$ the so-called glue code (a subgroup of the discriminant group of $\left\langle \sR \right\rangle$) (cf. \cite{bs96}, $4.4$, $4.5$, $4.7$).
Given a pair $(\sR,\CL)$, the associated lattice $L$ is constructed by $L = \left\langle \sR \right\rangle + \left\langle x \| \overline{x} \in \CL \right\rangle$.
We will refer to this result in Lemma $3.4$, Lemma $3.5$ and Lemma $3.9$.
\end{rem}

\begin{defn}Let $L$ be an integral lattice and $\CG$ its genus.
 \begin{itemize}
  \item[a)] We call $\CG$ \emph{totally-reflective} if each lattice in $\CG$ is reflective.
  \item[b)] The integral lattice $L$ is called \emph{totally-reflective} if its genus $\CG$ is totally-reflective.
 \end{itemize}
\end{defn}
One can deduce from the work of Biermann \cite{bie81}, that there are only finitely many totally-reflective genera in any fixed dimension (cf. \cite{sw92}, Theorem 1.4). 
Furthermore, Esselmann proved in \cite{ess96} that $20$ is the largest dimension of totally-reflective genera, thus a classification is possible (at least) in priniciple.
With the present paper we contribute to this problem by classifying the dimensions $3$ and $4$.
\section{Bounds for strongly square free, totally-reflective genera}
In this section we prove the results concerning step 1 and 2 of the general strategy. 
The basic idea is to compare the whole mass of a lattice $L$ with the part coming from the reflective lattices within $\CG(L)$. 
Since the latter quantity is crucial, we make the following
\begin{defn}
 Let $L$ be an integral lattice. We refer to
 \begin{align*}
  \mref(L) := \sum_{\substack{M \in \CG(L),\\\text{$M$ is reflective}}} \frac{1}{|O(M)|}
 \end{align*}
 as the \emph{reflective part} of the mass.
\end{defn}
An important (though trivial) observation is that $\mref(L) \leq \m(L)$ and $L$ is totally-reflective iff $\mref(L) = \m(L)$.
We will obtain our bounds by showing that the reflective part of the mass grows more slowly than the whole mass.
We begin with

\begin{lem}
  Let $D \in \NN$. Then
 \begin{itemize}
  \item[(a)] $\zeta_D(2) \geq \frac{\zeta(4)}{\zeta(2)} = \frac{\pi^2}{15}$.
  \item[(b)] $\zeta_{-D}(1) \leq 1+\frac{1}{2} \ln(D)$.
 \end{itemize}
 \begin{proof} 
  Part (a) follows from an elementary calculation which can be found in \cite{kl13}, $5.1$ and (b) is proved in \cite{wat79}, $5.10$.
 \end{proof}
\end{lem}

With part (a) of the previous lemma and the mass formula as stated in section 2.2, one can control the growth of $\m$ in the following sense.

\begin{lem}
 Let $L$ be a strongly square free lattice with determinant $d$.
 \begin{itemize}
  \item[(a)]For $\dim L = 3$ set
   \begin{align*}
    \M(L):=\frac{1}{6} \cdot \frac{1}{8} \cdot \prod_{\substack{p \mid d,\\p \neq 2}} \frac{p-1}{2}.
   \end{align*}
  \item[(b)]For $\dim L = 4$ set
   \begin{align*}
    \M(L):=\frac{1}{90}\cdot \frac{1}{24} \cdot \prod_{\substack{p \mid d,\\v_p(d) = 2 \\ p \neq 2}} \frac{1}{2}p^2\cdot\frac{p-1}{p+1} 
                       \cdot \prod_{\substack{p \mid d,\\v_p(d)=1 \\ p \neq 2}}\frac{1}{2}p^{\tfrac{3}{2}}.
   \end{align*}
 \end{itemize}
Then, in both cases, we have $\m(L) \geq \M(L)$.
\begin{proof}
 This follows directly from the mass formula and Lemma $3.2$.
\end{proof}
\end{lem}

It is somewhat more difficult to control the growth of $\mref$. For this, the following observation is helpful (recall Remark $2.4$).

\begin{lem}
 Let $L$ be an indecomposable reflective lattice with $\dim L \in \{2,3,4\}$. Then $O(L)$ only depends on the combinatorial class of $R(L)$.
 In particular, $O(L)$ does not depend on the glue-code nor the scaling. Referring to \cite{bs96}, $4.4$, $4.5$, $4.7$, we have
 \begin{itemize}
  \item[(a)] in dimension $2:$
  \renewcommand{\arraystretch}{1.2}
  \begin{center}
   $
   \begin{array}{c|c|c|c|c}
          & (a) & (b) & (c) & (d) \\
   \hline
   |O(L)| & 4 & 12 & 4 & 8
   \end{array}
   $
  \end{center}
  \item[(b)] in dimension $3:$
  \renewcommand{\arraystretch}{1.2}
  \begin{center}
   $
   \begin{array}{c|c|c|c|c|c}
          & (a) & (b) & (c) & (d) & (e)\\
   \hline
   |O(L)| & 8 & 8 & 16 & 48 & 48
   \end{array}
   $
  \end{center}
  \item[(c)] in dimension $4:$
  \renewcommand{\arraystretch}{1.2}
  \begin{center}
   $
   \begin{array}{c|c|c|c|c|c|c}
          & (a),(b),(c) & (d),(e) & (f),(g),(i) & (h) & (j),(k) & (l)\\
   \hline
   |O(L)| & 16 & 32 & 96 & 72 & 240 & 1152
   \end{array}
   $
  \end{center}
 \end{itemize}
 \begin{proof}
  Let $R(L)$ be decomposed as $R(L) = \xsR{\alpha_1}_1 \cdots \xsR{\alpha_k}_k$, with $\sR_i$ irreducible and let $O(R(L))$ be the stabilizer of $R(L)$ in $O(V)$.
  It is well known that $O(R(L)) \cong W(R(L)) \rtimes A(L)$, where $A(L)$ can be identified with a subgroup of the outer automorphism group of $R(L)$.
  It follows from \cite{bs96}, $4.4$, $4.5$, $4.7$, and the relation $W(R(L)) \subseteq O(L) \subseteq O(R(L))$ that $O(L) = W(R(L))$ holds in dimension $2$ and $3$
  since $A(L)$ is trivial, except for $\sA_2, \sA_3$ and $\sD_3$.
  Furthermore, $W(R(L))$ does not depend on the scaling since 
  $W(\xsR{\alpha_1}_1 \dots \xsR{\alpha_k}_k) \cong W(\xsR{\alpha_1}_1) \times \cdots \times W(\xsR{\alpha_k}_k) \cong W(\sR_1) \times \cdots \times W(\sR_k).$
  In dimension 4, one can use the same arguments, except in the cases $(h),(j)$ and $(k)$, where $A(L) \neq 1$. But there we have
  \begin{align*}
  A(L) = 
  \begin{cases}
   \ZZ/2\ZZ, & \text{in $(j)$ and $(k)$}, \\
   \ZZ/2\ZZ \times \ZZ/2\ZZ, & \text{in $(h),$} \\
  \end{cases}
  \end{align*}
  which obviously do not depend on scaling nor the glue-code.
  Weyl groups of unscaled, irreducible root systems, particularly their orders, are well known and can be found, for instance, in \cite{bou08}.
 \end{proof}
\end{lem}

\subsection{Bounds on the number of prime factors}
Regarding the statement on the number of prime factors outlined in step 1, the following way of controlling $\mref$ seems appropriate
(we will introduce an alternative way in the next subsection).
We recall that $\omega(d)$ (resp. $\Omega(d)$) refers to the number of (not necessarily) different prime factors of $d$.
\begin{lem}
 Let $L$ be a strongly square free lattice with determinant $d$.
 \begin{itemize}

  \item[a)] For $\dim L = 3$ we set
    \begin{align*}
     \Mref(L):=\sum_{\substack{x \mid d}}2^1 \cdot \frac{1}{2} \cdot \frac{2}{\pi} \cdot \prod_{p \mid x}\frac{1}{2}\sqrt{p} \cdot \left(1+\frac{1}{2}\ln(x)\right).
    \end{align*}

  \item[b)] For $\dim L = 4$ we set
    \begin{align*}
     \Mref(L) &:=3\frac{4^{\Omega(d)}}{16} + 2\frac{3^{\Omega(d)}}{32} + \frac{2^{\Omega(d)}}{72} + 3\frac{2^{\Omega(d)}}{96} + \frac{53}{5760} \\
                         &+ \sum_{\substack{x \mid d}} 2^{\omega(x)+1} \cdot \frac{1}{4} \cdot \frac{2}{\pi} 
                         \cdot  \left(1+\frac{1}{2}\ln(x)\right)
                         \cdot \prod_{\substack{p \mid x,\\v_p(d)=2}}\frac{2p}{2p-1} \cdot \prod_{\substack{p \mid x,\\v_p(d)=1}}\frac{1}{2}\sqrt{p},                     
    \end{align*}
 \end{itemize}
 Then, in both cases, we have $\mref(L) \leq \Mref(L)$.
 \begin{proof}
  $(b)$ We take a closer look on how the quantity $\mref$ is composed, more precisely we write
   \begin{align*}
    \mref(L) = \m^{(4)}_{\mathrm{ref}}(L) + \m^{(3)}_{\mathrm{ref}}(L) + \m^{(2)}_{\mathrm{ref}}(L),
   \end{align*}
   where $\m^{(4)}_{\mathrm{ref}}(L)$ refers to the part of the mass coming from the indecomposable, reflective lattices within the genus,
   $\m^{(3)}_{\mathrm{ref}}(L)$ refers to the part which comes from the reflective lattices with a $3$-dimensional, indecomposable component and
   $\m^{(2)}_{\mathrm{ref}}(L)$ means the contribution of the reflective lattices which have a $2$-dimensional, indecomposable compontent or which are diagonalisable.
   An estimate of $\m^{(4)}_{\mathrm{ref}}(L)$ is easily obtained by combining the classification theorem from \cite{bs96} and Lemma $3.4$.
   Considering that a bound for the number of possible isometry classes 
   for a given pair $(\sR,\CL)$ of determinant $d$ is $a^{\Omega(d)}$, where $a$ is the number of occurring scaling factors, we get
   \begin{align*}
     \m^{(4)}_{\mathrm{ref}}(L) \leq 3\frac{4^{\Omega(d)}}{16} + 2\frac{3^{\Omega(d)}}{32} + \frac{2^{\Omega(d)}}{72} + 3\frac{2^{\Omega(d)}}{96} + \frac{53}{5760}.
   \end{align*}
   As Berger showed in \cite{ber93}, a $4$-dimensional, strongly square free, reflective lattice can not have a $3$-dimensional, indecomposable component, so $\m^{(3)}_{\mathrm{ref}}(L) = 0$.
   For the estimate on $\m^{(2)}_{\mathrm{ref}}(L)$ we use the mass formula:
   \begin{align*}
    \m^{(2)}_{\mathrm{ref}}(L) &= \sum_{\substack{\text{$M \in \CG(L)$ refl,} \\M = M_1 \perp M_2,\\ \dim M_1=2}} \frac{1}{|O(M)|} 
                               \leq \frac{1}{4} \sum_{\substack{\text{$M \in \CG(L)$ refl,} \\M = M_1 \perp M_2,\\ \dim M_1=2}} \frac{1}{|O(M_1)|} \\                       
                               &\leq \sum_{\substack{x \mid d}} 2^{\omega(x)+1} \cdot \frac{1}{4} \cdot \frac{2}{\pi} 
                                \cdot  \left(1+\frac{1}{2}\ln(x)\right)
                                \cdot \prod_{\substack{p \mid d,\\v_p(d)=2}}\frac{2p}{2p-1} \cdot \prod_{\substack{p \mid d,\\v_p(d)=1}}\frac{1}{2}\sqrt{p},
   \end{align*}
   where $2^{\omega(x)+1}$ is an estimate for the number of $2$-dimensional, square free genera of determinant $x$. The estimate on $\m(M_1)$ follows from Proposition $2.2$ and Lemma $3.2. (b)$.

  $(a)$ Here we consider 
   \begin{align*}
    \mref(L) = \m^{(3)}_{\mathrm{ref}}(L) + \m^{(2)}_{\mathrm{ref}}(L).
   \end{align*}
  Again, $\m^{(3)}_{\mathrm{ref}}(L) = 0$, because $4$ divides the determinant of a $3$-dimensional, indecomposable, reflective lattice which is never the case for strongly square free lattices.
  For $\m^{(2)}_{\mathrm{ref}}(L)$ we get:
  \begin{align*}
   \m^{(2)}_{\mathrm{ref}}(L) &= \sum_{\substack{\text{$M \in \CG(L)$ refl,} \\M = M_1 \perp M_2,\\ \dim M_1=2}} \frac{1}{|O(M)|} 
                              \leq \frac{1}{2} \sum_{\substack{\text{$M \in \CG(L)$ refl,} \\M = M_1 \perp M_2,\\ \dim M_1=2}} \frac{1}{|O(M_1)|} \\                            
                              &\leq \sum_{\substack{x \mid d}} 2^1 \cdot \frac{1}{2} \cdot \frac{2}{\pi} \cdot \prod_{p \mid x}\frac{1}{2}\sqrt{p} \cdot \left(1+\frac{1}{2}\ln(x)\right).
  \end{align*}
  Notice that, unlike in the $4$-dimensional case, an estimate for the number of $2$-dimensional, square free genera of determinant $x$ is $2$.
  In the situation $\CG(M_1 \perp N_1) = \CG(M_2 \perp N_2)$ with $\dim M_i =2, \dim N_i=1$ and $\det N_i = x$ it follows from Witt's cancellation theorem for nondyadic local rings that
  $M_1 \otimes_{\ZZ} \ZZ_p \cong M_2 \otimes_{\ZZ} \ZZ_p$ for $p \neq 2$, whereas $M_1 \otimes_{\ZZ} \ZZ_2$ and $M_2 \otimes_{\ZZ} \ZZ_2$ can be distinguished only through the parity, 
  cf. \cite{oma00}, $92:3, 93:16$.
 \end{proof}
\end{lem}
By combining Lemma $3.3$ and Lemma $3.5$, we see that a strongly square free, totally-reflective lattice satisfies the condition $\Mref(L) / \M(L) \geq 1$.
Actually, the estimates on $\mref(L)$ and $\m(L)$ depend only on the determinant of $L$, so we may write $\M(L) = \M(\det L)$ and $\Mref(L) = \Mref(\det L)$.

To show that the number of prime factors can not be arbitrarily large, 
it is important to investigate how the ratio of $\Mref$ and $\M$ behaves when prime factors are appended to the determinant.
Notice that both $\Mref$ and $\M$ tend to $\infty$ when the number of prime factors increases.
\begin{lem} Let $d \in \mathbb{N}$ and $q \in \PP$.
 \begin{itemize}
  \item[(a)] In dimension $3$ we have 
   \begin{align*}
     \frac{\Mref(d)}{\M(d)} \geq \frac{\Mref(dq)}{\M(dq)}
   \end{align*}
  if $q \geq 17$ and $\Omega(d) \geq 2$.
  \item[(b)]In dimension $4$ we have 
   \begin{align*}
     \frac{\Mref(d)}{\M(d)} \geq \frac{\Mref(dq^2)}{\M(dq^2)}
   \end{align*}
  if $q \geq 7$ and $\Omega(d) \geq 2$, as well as
   \begin{align*}
    \frac{\Mref(d)}{\M(d)} \geq \frac{\Mref(dq)}{\M(dq)}
   \end{align*}
  if $q \geq 5$ and $\Omega(d) \geq 2$.
 \end{itemize}
 \begin{proof}
  We proof part (b) in detail so that the general idea will be clear. Part (a) is proven analogously.
  It is helpful to consider the quantity $\Mref$ piecewise. We write 
  $\Mref(d) = \M^{(4)}_{\mathrm{ref}}(d) + \M^{(2)}_{\mathrm{ref}}(d)$ with
  \begin{align*}
   \M^{(4)}_{\mathrm{ref}}(d):=3\frac{4^{\Omega(d)}}{16} + 2\frac{3^{\Omega(d)}}{32} + \frac{2^{\Omega(d)}}{72} + 3\frac{2^{\Omega(d)}}{96} + \frac{53}{5760}
  \end{align*}
  and
  \begin{align*}
   \M^{(2)}_{\mathrm{ref}}(d):=\sum_{\substack{x \mid d}} 2^{\omega(x)+1} \cdot \frac{1}{4} \cdot \frac{2}{\pi} 
                             \cdot  \left(1+\frac{1}{2}\ln(x)\right)
                             \cdot \prod_{\substack{p \mid d,\\v_p(d)=2}}\frac{2p}{2p-1} \cdot \prod_{\substack{p \mid d,\\v_p(d)=1}}\frac{1}{2}\sqrt{p}.
  \end{align*}
  The case $\M^{(4)}_{\mathrm{ref}}$ is easily done. From the mass formula we get the equivalence
  \begin{align*}
     \frac{\M^{(4)}_{\mathrm{ref}}(d)}{\M(d)} \geq \frac{\M^{(4)}_{\mathrm{ref}}(dq^2)}{\M(dq^2)} 
     \Longleftrightarrow \M^{(4)}_{\mathrm{ref}}(d) \cdot \frac{q^2(q-1)}{2 (q+1)} \geq \M^{(4)}_{\mathrm{ref}}(dq^2)
  \end{align*}
 and because of 
 \begin{align*}
  \M^{(4)}_{\mathrm{ref}}(dq^2) \geq 4^2 \cdot \M^{(4)}_{\mathrm{ref}}(d)
 \end{align*}
 the inequality on the right hand side is satisfied when $\tfrac{q^2(q-1)}{2 (q+1)} \geq 4^2$,
 that is when $q \geq 7$.
 
 The case $\M^{(2)}_{\mathrm{ref}}$ is more difficult. Let $p_1$ be an arbitrary prime dividing $d$. 
 Define $D(d):=\{ x \in \NN \| \text{$d$ is divisible by $x$}\}$.
 For the first part of (b) the starting point is the decomposition
 \begin{align}
  D(dq^2) & =  D(d)\setminus \{p_1\} \\
                & \cupdot  q \cdot \left(D(d)\setminus \{p_1,1\} \right) \\
                & \cupdot q^2 \cdot \left(D(d)\setminus \{p_1,1\} \right) \\
                & \cupdot \left\{q, q^2, p_1, q p_1, q^2  p_1 \right\}.
 \end{align}
 The idea is to compare the summands of $\M^{(2)}_{\mathrm{ref}}(dq^2)$ corresponding to the right hand side of $(2), (3)$ and $(4)$ (resp.
 $(5)$) with the summands of $\M^{(2)}_{\mathrm{ref}}(d)$ corresponding to $D(d)\setminus \{p_1\}$ (resp. $p_1$).
 For all $x \in D(d) \setminus \{p_1\}$ the mass formula implies
 \begin{align*}
  & \left( 2^{\omega(x)+1} \cdot \frac{1}{4} \cdot \frac{2}{\pi} 
    \cdot  \left(1+\frac{1}{2}\ln(x)\right)
    \cdot \prod_{\substack{p \mid d,\\v_p(d)=2}}\frac{2p}{2p-1} \cdot \prod_{\substack{p \mid d,\\v_p(d)=1}}\frac{1}{2}\sqrt{p}\right)
    \cdot \frac{q^2(q-1)}{2 (q+1)} \\
  & \geq  2^{\omega(x)+1} \cdot \frac{1}{4} \cdot \frac{2}{\pi} 
           \cdot  \left(1+\frac{1}{2}\ln(x)\right)
           \cdot \prod_{\substack{p \mid d,\\v_p(d)=2}}\frac{2p}{2p-1} \cdot \prod_{\substack{p \mid d,\\v_p(d)=1}}\frac{1}{2}\sqrt{p} \\          
  & + 2^{\omega(qx)+1} \cdot \frac{1}{4} \cdot \frac{2}{\pi} 
      \cdot  \left(1+\frac{1}{2}\ln(qx)\right)
      \cdot \prod_{\substack{p \mid d,\\v_p(d)=2}}\frac{2p}{2p-1} \cdot \prod_{\substack{p \mid d,\\v_p(d)=1}}\frac{1}{2}\sqrt{p} \cdot \frac{1}{2}\sqrt{q}\\
  & + 2^{\omega(q^2x)+1} \cdot \frac{1}{4} \cdot \frac{2}{\pi} 
      \cdot  \left(1+\frac{1}{2}\ln(q^2x)\right)
      \cdot \prod_{\substack{p \mid d,\\v_p(d)=2}}\frac{2p}{2p-1} \cdot \frac{2q}{2q-1} \cdot \prod_{\substack{p \mid d,\\v_p(d)=1}}\frac{1}{2}\sqrt{p},
 \end{align*}
 which is equivalent to 
 \begin{align*}
  & \left(\frac{q^2(q-1)}{2 (q+1)}-1\right) \cdot \left(1+\frac{1}{2}\ln(x)\right) \\
  & \geq \left( \sqrt{q}+\frac{\sqrt{q}}{2} \ln(q) + \frac{4q}{2q-1}+\frac{2q}{2q-1} \ln(q^2) \right) + \left( \frac{\sqrt{q}}{2} + \frac{2q}{2q-1} \right) \ln(x).
 \end{align*}
 Since $x \geq 1$, and thus $\ln(x) = 0$ or $\ln(x) \geq 1$, the latter inequality is true for $q \ge 7$. For $x = p_1$ we obtain
  \begin{align*}
    & \left( 2^{\omega(p_1)+1} \cdot \frac{1}{4} \cdot \frac{2}{\pi} 
       \cdot  \left(1+\frac{1}{2}\ln(p_1)\right)
       \cdot \frac{1}{2}\sqrt{p_1} \right) \cdot \frac{q^2(q-1)}{2 (q+1)}  \\
   & + 2^{\omega(p_1)+1} \cdot \frac{1}{4} \cdot \frac{2}{\pi} 
       \cdot  \left(1+\frac{1}{2}\ln(p_1)\right)
       \cdot \frac{1}{2}\sqrt{p_1} \\
   & +  2^{\omega(q)+1} \cdot \frac{1}{4} \cdot \frac{2}{\pi} 
       \cdot  \left(1+\frac{1}{2}\ln(q)\right)
       \cdot \frac{1}{2}\sqrt{q} \\
   & +  2^{\omega(q^2)+1} \cdot \frac{1}{4} \cdot \frac{2}{\pi} 
       \cdot  \left(1+\frac{1}{2}\ln(q^2)\right)
       \cdot \frac{2q}{2q-1} \\
   & +  2^{\omega(qp_1)+1} \cdot \frac{1}{4} \cdot \frac{2}{\pi} 
       \cdot  \left(1+\frac{1}{2}\ln(qp_1)\right)
       \cdot \frac{1}{2}\sqrt{q} \cdot \frac{1}{2} \sqrt{p_1} \\
   & +  2^{\omega(q^2p_1)+1} \cdot \frac{1}{4} \cdot \frac{2}{\pi} 
       \cdot  \left(1+\frac{1}{2}\ln(q^2p_1)\right)
       \cdot  \frac{2q}{2q-1} \cdot \frac{1}{2}\sqrt{p_1}.
  \end{align*}
 It is easy to check that this is true for all $p_1 \geq 2$ and $q \geq 7$.
 The second part of (b) works similarly, but uses the decompostion
 \begin{align*}
    D(dq) & =  D(d) \setminus \{p_1\} \\
          & \cupdot  q \cdot \left(D(d)\setminus \{p_1\}\right) \\
          & \cupdot \left\{q, p_1, q p_1 \right\}.
 \end{align*}
 \end{proof}
\end{lem}
The second lemma of this section clarifies the behavior of $\Mref / \M$ if the number of prime factors is fixed while the prime numbers increase.
Again, a priori that is not clear since $\Mref$ and $\M$ tend to $\infty$ when the primes increase.
\begin{lem}
 In both dimensions $\frac{\Mref(d)}{\M(d)}$ is monotonically decreasing in each prime factor of $d$.
 \begin{proof}
  Let $q$ be a prime with $v_q(d) = 2$ (the easier case $v_p(d) = 1$ is proven analogously). 
  Like in the previous proof we consider $\M^{(4)}_{\mathrm{ref}}(d)$ and $\M^{(2)}_{\mathrm{ref}}(d)$ separated.
  Since $\M^{(4)}_{\mathrm{ref}}(d)$ only depends on the number of prime factors (and not on the primes itself), 
  it is clear that $\frac{\M^{(4)}_{\mathrm{ref}}(d)}{\M(d)}$ decreases monotonically. In order to deal with $\M^{(2)}_{\mathrm{ref}}(d)$, we have to isolate the effect of $q$.
  For this the following decomposition is helpful:
  \begin{align*}
   D(d) = D(\tfrac{d}{q^2}) \cupdot q D(\tfrac{d}{q^2}) \cupdot q^2 D(\tfrac{d}{q^2}).
  \end{align*}
  Thus we can write
  \begin{align*}
  & \M^{(2)}_{\mathrm{ref}}(d)  = \sum_{\substack{x \mid d}} 2^{\omega(x)+1} \cdot \frac{1}{4} \cdot \frac{2}{\pi} 
     \cdot  \left(1+\frac{1}{2}\ln(x)\right)
     \cdot \prod_{\substack{p \mid d,\\v_p(d)=2}}\frac{2p}{2p-1} \cdot \prod_{\substack{p \mid d,\\v_p(d)=1}}\frac{1}{2}\sqrt{p} \\
   & =   \underbrace{\sum_{\substack{x \mid d/q}} 2^{\omega(x)+1} \cdot \frac{1}{4} \cdot \frac{2}{\pi} 
         \cdot  \left(1+\frac{1}{2}\ln(x)\right)
         \cdot \prod_{\substack{p \mid d/q,\\v_p(d/q)=2}}\frac{2p}{2p-1} \cdot \prod_{\substack{p \mid d/q,\\v_p(d/q)=1}}\frac{1}{2}\sqrt{p}}_{\text{$:=c_1$, does not depend on $q$}} \\
   & +   \sum_{x \mid d/q} \underbrace{2^{\omega(x)+1} \cdot \frac{1}{4} \cdot \frac{2}{\pi} 
         \cdot \prod_{\substack{p \mid d/q,\\v_p(d/q)=2}}\frac{2p}{2p-1} \cdot \prod_{\substack{p \mid d/q,\\v_p(d/q)=1}}\frac{1}{2}\sqrt{p} \cdot \frac{1}{2}}_{\text{$:=c_2(x)$, does not depend on $q$}}
         \cdot \sqrt{q} \cdot \left(1+\frac{1}{2}\ln(qx)\right)\\
   & +   \sum_{x \mid d/q} \underbrace{2^{\omega(x)+1} \cdot \frac{1}{4} \cdot \frac{2}{\pi} 
          \cdot \prod_{\substack{p \mid d/q,\\v_p(d/q)=2}}\frac{2p}{2p-1} \cdot \prod_{\substack{p \mid d/q,\\v_p(d/q)=1}}\frac{1}{2}\sqrt{p}}_{\text{$:=c_3(x)$, does not depend on $q$}}
          \cdot \frac{2q}{2q-1} \cdot \left(1+\frac{1}{2}\ln(q^2x)\right).
  \end{align*}
  Furthermore, by isolating the effect of $q$ in $\M(d)$,
  \begin{align*}
      \M(d) = \underbrace{\frac{1}{90}\cdot \frac{1}{24} \cdot \prod_{\substack{p \mid d/q,\\v_p(d) = 2}} \frac{1}{2}p^2\cdot\frac{p-1}{p+1} 
                   \cdot \prod_{\substack{p \mid d,\\v_p(d)=1}}\frac{1}{2}p^{\tfrac{3}{2}} \cdot \frac{1}{2}}_{\text{$:=c$, does not depend on $q$}} 
                   \cdot q^2\cdot\frac{q-1}{q+1}
  \end{align*}
  we can treat 
  \begin{align*}
   \frac{\M^{(2)}_{\mathrm{ref}}(d)}{\M(d)} = \frac{c_1 + \sum_{x} c_2(x)\cdot \sqrt{q} \cdot \left(1+\tfrac{1}{2}\ln(qx)\right) + 
                                              \sum_{x} c_3(x) \cdot \tfrac{2q}{2q-1} \cdot \left(1+\tfrac{1}{2}\ln(q^2x)\right)}{c \cdot q^2\cdot\frac{q-1}{q+1}}
  \end{align*}
  as a differentiable function of $q$. Then one easily checks that the first derivative after $q$ is $<0$.
 \end{proof}
\end{lem}
The main theorem of this subsection is now a direct consequence of Lemma $3.6$ and Lemma $3.7$. 
\begin{thm}
 Let $L$ be a strongly square free, totally-reflective lattice.
\begin{itemize}
 \item[(a)] Let $\dim L=3$ and $\det L =q_1 \cdots q_s$. Then $s \leq 10$.
 \item[(b)] Let $\dim L=4$ and $\det L =p_1^2 \cdots p_r^2 q_1 \cdots q_s$. Then $r \leq 9$ and $s \leq 8-r$.
\end{itemize}
 \begin{proof}
 We have to decide when the necessary condition $\frac{\Mref}{\M} \geq 1 $ is violated.

 (a) Let $L$ be a strongly square free, totally-reflective lattice with $\det L = q_1 \cdots q_s$ and $s \geq 11$. Assume the prime factors are ordered such that $q_1 < \dots < q_s$.
 We start with the observation
 \begin{align*}
  \frac{\Mref(2 \cdot 3 \cdot 5 \cdot 7 \cdot 11 \cdot 13 \cdot 17 \cdot 19 \cdot 23 \cdot 29 \cdot 31)}{\M(2 \cdot 3 \cdot 5 \cdot 7 \cdot 11 \cdot 13 \cdot 17 \cdot 19 \cdot 23 \cdot 29 \cdot 31)} < 1.
 \end{align*}
 Using the monotony statement of Lemma $3.7$, we get
 \begin{align*}
  \frac{\Mref(q_1 \cdot q_2 \cdot q_3 \cdot q_4 \cdot q_5 \cdot q_6 \cdot q_7 \cdot q_8 \cdot q_9 \cdot q_{10} \cdot q_{11})}
       {\M(q_1 \cdot q_2 \cdot q_3 \cdot q_4 \cdot q_5 \cdot q_6 \cdot q_7 \cdot q_8 \cdot q_9\cdot q_{10} \cdot q_{11})} < 1,
 \end{align*}
 which is possible since $s \geq 11$. Now we apply Lemma $3.6$ (a) and see that
  \begin{align*}
  \frac{\Mref(q_1 \cdots q_{11} \cdot q_{12} \cdots q_s)}{\M(q_1 \cdots q_{11} \cdot q_{12} \cdots q_s)} = \frac{\Mref(\det L)}{\M(\det L)} < 1.
 \end{align*}
 Thus $L$ is not totally-reflective.

 (b) Let $L$ be a strongly square free, totally-reflective lattice of determinant $\det L = p_1^2 \cdots p_r^2 q_1 \cdots q_s$ with $r \geq 10$. 
 First we prove the statement regarding $r$.
 Assume $p_1 < \dots < p_r$ and $q_1 < \dots < q_s$. We have 
 \begin{align*}
  \frac{\Mref(2^2 \cdot 3^2 \cdot 5^2 \cdot 7^2 \cdot 11^2 \cdot 13^2 \cdot 17^2 \cdot 19^2 \cdot 23^2 \cdot 29^2)}
       {\M(2^2 \cdot 3^2 \cdot 5^2 \cdot 7^2 \cdot 11^2 \cdot 13^2 \cdot 17^2 \cdot 19^2 \cdot 23^2 \cdot 29^2)} < 1
 \end{align*}
 and with Lemma $3.7$
 \begin{align*}
  \frac{\Mref(p_1^2 \cdot p_2^2 \cdot p_3^2 \cdot p_4^2 \cdot p_5^2 \cdot p_6^2 \cdot p_7^2 \cdot p_8^2 \cdot p_9^2 \cdot p_{10}^2)}
       {\M(p_1^2 \cdot p_2^2 \cdot p_3^2 \cdot p_4^2 \cdot p_5^2 \cdot p_6^2 \cdot p_7^2 \cdot p_8^2 \cdot p_9^2 \cdot p_{10}^2)} < 1.
 \end{align*}
 Applying Lemma $3.6$ part (a) and (b), we get
 \begin{align*}
  \frac{\Mref(p_1^2 \cdots p_{10}^2 \cdot p_{11}^2 \cdots p_{r}^2 \cdot q_1 \cdots q_s)}
       {\M(p_1^2 \cdots p_{10}^2 \cdot p_{11}^2 \cdots p_{r}^2 \cdot q_1 \cdots q_s)} = \frac{\Mref(\det L)}{\M(\det L)} < 1.
 \end{align*}
 Thus $L$ is not totally-reflective.

 To prove the statement concerning $s$, we fix the number of quadratic prime factors $r \leq 9$. Let $\det L = p_1^2 \cdots p_r^2 q_1 \cdots q_s$ with $s \geq 9-r$.
 Define $\PP(r+s)$ to be the finite set consisting of the first $r+s$ primes. For each combination $(\tilde{p}_1, \dots, \tilde{p}_r,\tilde{q}_1, \dots, \tilde{q}_s) \in \PP(r+s)^{r+s}$ 
 with $\tilde{p}_1 < \dots < \tilde{p}_r$ and $\tilde{q}_1 < \dots < \tilde{q}_s$ we have
  \begin{align*}
  \frac{\Mref(\tilde{p}_1^2 \cdots \tilde{p}_r^2 \tilde{q}_1 \cdots \tilde{p}_s)}
       {\M(\tilde{p}_1^2 \cdots \tilde{p}_r^2 \tilde{q}_1 \cdots \tilde{p}_s)} < 1.
 \end{align*}
 Then Lemma $3.6$ and Lemma $3.7$ implie
   \begin{align*}
  \frac{\Mref(p_1^2 \cdots p_r^2 q_1 \cdots q_s)}
       {\M(p_1^2 \cdots p_r^2 q_1 \cdots q_s)} < 1.
 \end{align*}
 \end{proof}
\end{thm}

\subsection{Bounds on the prime factors}
In the previous subsection we used the mass formula to estimate the part of the mass coming from those reflective lattices that decompose into one- and two-dimensional sublattices.
For the purpose of this subsection (step $2$ of the general strategy) a different approach seems more appropriate.
\begin{lem}
 Let $L$ be a strongly square free lattice with determinant $d$.
 \begin{itemize}

  \item[a)] For $\dim L = 3$ we set
    \begin{align*}
     \Nref(L):=\sum_{x \mid d} \frac{1}{2} \left( 2\frac{2^{\Omega(x)}}{4} + \frac{5}{24} \right).
    \end{align*}

  \item[b)] For $\dim L = 4$ we set
    \begin{align*}
     \Nref(L) &:=3\frac{4^{\Omega(d)}}{16} + 2\frac{3^{\Omega(d)}}{32} + \frac{2^{\Omega(d)}}{72} + 3\frac{2^{\Omega(d)}}{96} + \frac{53}{5760} \\
                         &+ \sum_{x \mid d} \frac{1}{4} \left( 2\frac{2^{\Omega(x)}}{4} + \frac{5}{24} \right)\cdot \left( 2\frac{2^{\Omega(d/x)}}{4} + \frac{5}{24} \right).
    \end{align*}
 \end{itemize}
 Then, in both cases, we have $\mref(L) \leq \Nref(L)$.
 \begin{proof}
  (b) Like in the proof of Lemma $3.5$ we start with the equation
  \begin{align*}
     \mref(L) = \m^{(4)}_{\mathrm{ref}}(L) + \m^{(3)}_{\mathrm{ref}}(L) + \m^{(2)}_{\mathrm{ref}}(L).
  \end{align*}
  The estimates on $\m^{(3)}_{\mathrm{ref}}(L)$ and $\m^{(4)}_{\mathrm{ref}}(L)$ are the same as in $3.5$. But this time we rather use the classification theorem from \cite{bs96} than the mass formula
   to control $\m^{(2)}_{\mathrm{ref}}(L)$. Consider a reflective lattice $M = M_1 \perp M_2$ with $\det M_1 =x$ and the associated pair $(\sR(M_1),\CL(M_1))$.
   The classification \cite{bs96}, $4.4$, $4.5$, $4.7$, implies that a bound for the number of possible isometry classes for $M_1$ is $a^{\Omega(x)}$, 
   where $a$ is the number of occurring scaling factors. Furthermore, each of the $a^{\Omega(x)}$ isometry classes can occur $a^{\Omega(d/x)}$ times, thus
   \begin{align*}
    \m^{(2)}_{\mathrm{ref}}(L) & = \sum_{\substack{\text{$M \in \CG(L)$ refl,} \\M = M_1 \perp M_2,\\ \dim M_1=2}} \frac{1}{|O(M)|} 
                               \leq \frac{1}{4} \sum_{\substack{\text{$M \in \CG(L)$ refl,} \\M = M_1 \perp M_2,\\ \dim M_1=2}} \frac{1}{|O(M_1)|} \\                       
                               & \leq \sum_{x \mid d} \frac{1}{4} \left( 2\frac{2^{\Omega(x)}}{4} + \frac{5}{24} \right)\cdot \left( 2\frac{2^{\Omega(d/x)}}{4} + \frac{5}{24} \right).
   \end{align*}

  (a) When considering a $3$-dimensional lattice $M$ decomposed as $M = M_1 \perp M_2$, such that $\dim M_1 =2, \dim M_2 =1$, the estimate for the number of possible isometry classes for $M_1$
      works similarly. Each isometry class, however, occurs only one time which follows from Witt's cancelation theorem. Thus 
  \begin{align*} \m^{(2)}_{\mathrm{ref}}(L) &= \sum_{\substack{\text{$M \in \CG(L)$ refl,} \\M = M_1 \perp M_2,\\ \dim M_1=2}} \frac{1}{|O(M)|} 
                                             \leq \frac{1}{2} \sum_{\substack{\text{$M \in \CG(L)$ refl,} \\M = M_1 \perp M_2,\\ \dim M_1=2}} \frac{1}{|O(M_1)|} \\
                                            & \leq \sum_{x \mid d} \frac{1}{2} \left( 2\frac{2^{\Omega(x)}}{4} + \frac{5}{24} \right).
  \end{align*}
 \end{proof}
\end{lem}
In analogy to $\Mref$, the quantity $\Nref$ only depends on the prime factors of the determinant, thus the formulation in the following lemma makes sense.
\begin{lem}
 In both dimensions $\frac{\Nref(d)}{\M(d)}$ is monotonically decreasing in each prime factor of $d$.
 \begin{proof}
   The enumerator $\Nref(d)$ only depends on the number of prime factors of $d$ and, as shown in the proof of Lemma $3.7$, the denominator $\M(d)$ is monotonically increasing in each prime factor.
 \end{proof}
\end{lem}
As a first application we can slightly improve our bounds for the number of prime factors in dimension $3$.
\begin{cor}
 Let $L$ be a strongly square free, totally-reflective lattice with $\dim L=3$ and $\det L =q_1 \cdots q_s$. Then $s \leq 9$.
 \begin{proof}
  This follows from 
  \begin{align*}
   \frac{\Nref(2 \cdot 3 \cdot 5 \cdot 7 \cdot 11 \cdot 13 \cdot 17 \cdot 19 \cdot 23 \cdot 29)}{\M(2 \cdot 3 \cdot 5 \cdot 7 \cdot 11 \cdot 13 \cdot 17 \cdot 19 \cdot 23 \cdot 29)} < 1
  \end{align*}
  and Lemma $3.10$.
 \end{proof}
\end{cor}

Keeping in mind that a strongly square free, totally-refelctive lattice $L$ fulfills the condition $\Nref(L) / \M(L) \geq 1$, 
Lemma $3.10$ provides (from the computational point of view) strong bounds on the prime factors of $\det L$. 
\begin{thm}
 Let $L$ be a strongly square free, totally-reflective lattice.
\begin{itemize}
 \item[(a)] Let $\dim L=3$ and $\det L =q_1 \cdots q_s$ with $s \leq 9$. Assume $q_1 < \dots < q_s$. Then
 \medskip
 \begin{center}
  $
 \begin{array}{c|c|c|c|c|c|c|c|c|c|}
  & q_1 & q_2 & q_3 & q_4 & q_5 & q_6 & q_7 & q_8 & q_9\\
  \hline
  \leq & 89 & 257 & 733 & 1063 & 1033 & 607 & 293 & 113 & 37
 \end{array}.
  $
 \end{center}
 \medskip
 \item[(b)] Let $\dim L=4$ and $\det L =p_1^2 \cdots p_r^2 q_1 \cdots q_s$ with $r \leq 9$ and $s \leq 8-r$. Assume $p_1 < \dots < p_r$ and $q_1 < \dots < q_s$. Then
 \medskip
 \begin{center}
  $
 \begin{array}{c|c|c|c|c|c|c|c|c|c|}
  & p_1 & p_2 & p_3 & p_4 & p_5 & p_6 & p_7 & p_8 & p_9 \\
  \hline
  \leq & 191 & 661 & 1601 & 2069 & 1831 & 997 & 449 & 157 & 47
 \end{array},
  $
 \end{center}
 \medskip
 \begin{center}
  $
 \begin{array}{c|c|c|c|c|c|c|c|c|}
  & q_1 & q_2 & q_3 & q_4 & q_5 & q_6 & q_7 & q_8\\
  \hline
  \leq & 11287 & 6427 & 3613 & 1597 & 653 & 229 & 67 & 19
 \end{array}.
  $
 \end{center}
\end{itemize}
 \begin{proof}
  By using Lemma $3.10$ we can repeatedly increase a prime factor (and thus decrease the function $\Nref/\M$) 
  until the necessary condition $\Nref(\det L) / \M(\det L) \geq 1$ is violated.
 \end{proof}
\end{thm}

\begin{rem}\ 
 \begin{itemize}
  \item[(1)] With the help of Theorem $3.8$ and Theorem $3.12$, the enumeration of all strongly square free, totally-reflective genera can be carried out computationally.
       One produces all genera up to the given bounds and checks whether they are totally-reflective. The number of possible genera is finite since we deal with strongly square free lattices.
       The list is included in section $6$.
  \item[(2)] It turns out in section $6$ that the largest occuring value for the number of prime factors is 
             \begin{align*}
              (r,s)=
               \begin{cases}
                (3,3), & \text{in dimension $4$}, \\
                (0,4), & \text{in dimension $3$}. 
               \end{cases}
              \end{align*}
             The largest prime factor $p$ occuring in dimension $4$ is
             \begin{align*}
              p=
               \begin{cases}
                13, & \text{if $v_p(\det) = 2$}, \\
                17, & \text{if $v_p(\det) = 1$}. 
               \end{cases}
              \end{align*}
              and $p = 23$ in dimension 3.
 \end{itemize}
\end{rem}

\section{Completing the classification}
Now, after we found all strongly square free, totally-reflective genera, our next goal is to gradually weaken the restriction ``strongly square free''.
This is done in two steps. First, we drop the assumption ``strong'' by applying the partial dualization operator (this is easy). 
In a second step we drop the assumption ``square free'' which turnes out to be somewhat more difficult. 
We need to clarify for which primes we have to consider pre-images under the Watson transformation.

\subsection{From strongly square free to square free}
\begin{defn}
Let $L$ be an integral lattice and $p \in \PP$. The \emph{partial dual} of $L$ at $p$ is defined as $\D_p(L):={}^{p}\!{(\tfrac{1}{p}L \cap L^{\#}})$.
\end{defn}
In contrast to the usual dual operator, the partial dual operator only dualizes the lattice at the prime spot $p$. That means
\begin{align*}
 \D_p(L) \otimes_{\ZZ} \ZZ_q =
 \begin{cases}
  {}^{p}\!{\left(L_q^{\#}\right)}, & \text{if $q = p$}, \\
  {}^{p}\!{L_q},      & \text{if $q \neq p$}.
 \end{cases}
\end{align*}
This has the following effect on the Jordan decomposition of a square free lattice $L \otimes_{\ZZ} \ZZ_p = L_0 \perp {}^{p}{L_1}$:
\begin{align*}
 \D_p(L) \otimes_{\ZZ} \ZZ_p = L_1 \perp {}^{p}{L_0}.
\end{align*}
Thus, starting with a strongly square free lattice, one can construct a (not necessarily strongly) square free, primitive lattice by applying $\D_p$ for $p \mid \det L$ (and vice versa).
For a set of primes $I:=\left\{p_1,\cdots,p_k \right\} \subseteq \PP$ we use the abbreviation $\D_I:=\D_{p_1} \circ \cdots \circ \D_{p_k}$ (where $\D_\emptyset := \mathrm{id}$)
which is well-defined since two partial dual operators with respect to different primes commute.
Clearly, $\D_p$ extends to a well-defined bijective function $\CG(L) \longrightarrow \D_p(\CG(L))=\CG(\D_p(L))$.

The next lemma shows that the partial dual behaves ``well'' relative to the property ``totally-reflective''.
\begin{lem}
 Let $L$ be an integral lattice. Then $L$ is totally-reflective if and only if $\D_p(L)$ is totally-reflective.
 \begin{proof}
  Recall that the property ``reflective'' can be characterized by the action of $W(L)$ on $V$. Since $W(L)=W(\D_p(L))$, it is clear that $L$ is reflective iff $\D_p(L)$ is reflective.
  Thus, the assertion follows from the bijectivity of $\D_p: \CG(L) \longrightarrow \D_p(\CG(L))$.
 \end{proof}
\end{lem}
\begin{thm}
 Let $\CT_{n}$ be the set of all strongly square free, totally-reflective genera in dimension $n \in \{3,4\}$. Let $\CP(d)$ be the power set of the set of all prime factors of $d:=\det \CG$.
 Then 
 \begin{align*}
    \bigcup_{\CG \in \CT_n} \bigcup_{I \in \CP(d)} \left\{\D_I(\CG)\right\} 
 \end{align*}
 is the set of all square free, totally-reflective, primitve genera in dimension $n$.
\begin{proof}
 This is a consequence of Lemma $4.2$, the bijectivity of $\CG(L) \longrightarrow \D_p(\CG(L))$ and the above discussion.
\end{proof}
\end{thm}

\subsection{From square free to all}
The techniques we will use in this subsection are based on the following definition, going back to Watson, cf. \cite{wat62}, \cite{wat73}.
\begin{defn}
Let $L$ be an integral lattice and $p \in \PP$. The \emph{Watson transformation} of $L$ at $p$ is defined as $\E_p(L):=L + (\tfrac{1}{p}L \cap pL^{\#})$.
\end{defn}
The usefulness of $\E_p$ becomes clear when we consider its effect on the Jordan decomposition. Let $L$ be an integral lattice with
$L_p = L_0 \perp {}^{p}\!{L_1} \perp \dots \perp {}^{p^r}\!{L_r}$. Then
\begin{align*}
 \E_p(L) \otimes_{\ZZ} \ZZ_q =
 \begin{cases}
  \left(L_0 \perp L_2\right) \perp {}^{p}\!{\left(L_1 \perp L_3\right)} \perp {}^{p^2}\!{L_4} \perp \dots \perp {}^{p^{r-2}}\!{L_r}  , & \text{if $q = p$}, \\
  L \otimes_{\ZZ} \ZZ_q,      & \text{if $q \neq p$}.
 \end{cases}
\end{align*}
Hence, after repeatedly applying the Watson transformation, a primitive lattice transforms into a square free, primitive lattice. 
Similar to the partial dual,
$\E_p$ extends to a well-defined surjective function $\CG(L) \longrightarrow \E_p(\CG(L)) = \CG(\E_p(L))$.
\begin{lem}
 Let $L$ be a totally-reflective lattice. Then $\E_p(L)$ is totally-reflective.
 \begin{proof}
 The assertion implies that $W(L)$ has no nonzero fixed vectors, thus neither has $W(\E_p(L))$ since $W(L) \subseteq W(\E_p(L))$. Hence the assertion follows from the surjectivity of $\E_p$.
 \end{proof}
\end{lem}
It may happen that prime factors disappear from the determinant after applying the Watson transformation.
Thus, when calculating pre-images under $\E_p$, one has to decide which primes $p$ to consider (besides the prime factors of the determinant).
An answer to this question is given by the follwing two lemmata.
\begin{lem} Let $L$ be an integral lattice, $p$ an odd prime with $p \nmid \det L$ and $K \in \E_p(L)^{-1}$.
 \begin{itemize}
  \item[(a)] If $\dim L = 3$ then,       
             \begin{align*}
              \m(K) \geq \left(\left(\frac{1}{1+p^{-1}}\right)^2 p^2 \left(1-p^{-2}\right)\right) \cdot \m(L).
             \end{align*}             
  \item[(b)] If $\dim L = 4$ then,
             \begin{align*}
              \m(K) \geq \left( \frac{\zeta(4)}{2\zeta(2)^2} \left(\frac{1}{1+p^{-1}}\right)^2 p^3 \left(1-p^{-2}\right)\right) \cdot \m(L).
             \end{align*}
 \end{itemize}
 \begin{proof}
  This follows from the mass formula and Lemma $3.2$.
 \end{proof}
\end{lem}
Since we no longer deal only with strongly square free lattices, the following extended definition of $\Nref(L)$ becomes necessary:
\begin{align*}
 \Nref(L) := \N^{(4)}_{\mathrm{ref}}(L) + \N^{(3)}_{\mathrm{ref}}(L) + \N^{(2)}_{\mathrm{ref}}(L),
\end{align*}
with
\begin{align*}
 \N^{(4)}_{\mathrm{ref}}(L) &:=3\frac{4^{\Omega(d)}}{16} + 2\frac{3^{\Omega(d)}}{32} + \frac{2^{\Omega(d)}}{72} + 3\frac{2^{\Omega(d)}}{96} + \frac{53}{5760},\\
 \N^{(3)}_{\mathrm{ref}}(L) &:=\sum_{x \mid d}2\frac{3^{\Omega(x)}}{8} + \frac{2^{\Omega(x)}}{16}+\frac{1}{24},\\
 \N^{(2)}_{\mathrm{ref}}(L) &:=
 \begin{cases}
  \sum_{x \mid d} \frac{1}{4} \left( 2\frac{2^{\Omega(x)}}{4} + \frac{5}{24} \right)\cdot \left( 2\frac{2^{\Omega(d/x)}}{4} + \frac{5}{24} \right), & \text{if $\dim L =4$},\\
  \sum_{x \mid d} \frac{1}{2} \left( 2\frac{2^{\Omega(x)}}{4} + \frac{5}{24} \right), & \text{if $\dim L =3$}.
 \end{cases}
\end{align*}
It follows from a calculation similar to Lemma $3.9$ that $\mref(L) \leq \Nref(L)$, for a not necessarily strongly square free lattice $L$.
\begin{lem} Let $L$ be an integral lattice, $p$ an odd prime with $p \nmid \det L$ and $K \in \E_p(L)^{-1}$.
 \begin{itemize}
  \item[(a)] If $\dim L = 3$ then,       
             \begin{align*}
              \Nref(K) \leq 81 \cdot \Nref(L).
             \end{align*}             
  \item[(b)] If $\dim L = 4$ then,
             \begin{align*}
              \Nref(K) \leq 5103 \cdot \Nref(L).
             \end{align*}
 \end{itemize}
 \begin{proof}
  (b) Let $d:=\det L$. The assumption implies $K \otimes_{\ZZ} \ZZ_p = K_0 \perp {}^{p^{2}}\!{K_{2}}$, in particular $\det K = \det L \cdot p^{2\dim K_2}$, where $n_2:=\dim K_2 \leq 3$. Thus
  \begin{align*}
    \N^{(4)}_{\mathrm{ref}}(K) &=3\frac{4^{\Omega(dp^{2n_2})}}{16} + 2\frac{3^{\Omega(dp^{2n_2})}}{32} + \frac{2^{\Omega(dp^{2n_2})}}{72} 
                                              + 3\frac{2^{\Omega(dp^{2n_2})}}{96} + \frac{53}{5760} \\
                             & \leq 4^{2n_2}  \cdot \N^{(4)}_{\mathrm{ref}}(L).
  \end{align*}
 For the other two cases we consider the decomposition 
 \begin{align*}
  D(\det K) = \bigcupdot_{i=0}^{2n_2}p^i D(d).
 \end{align*}
 Hence
 \begin{align*}
   \N^{(3)}_{\mathrm{ref}}(K) &=\sum_{x \mid D(d)}2\frac{3^{\Omega(x)}}{8} + \frac{2^{\Omega(x)}}{16}+\frac{1}{24} \\
                                          &+\sum_{x \mid D(d)}2\frac{3^{\Omega(px)}}{8} + \frac{2^{\Omega(px)}}{16}+\frac{1}{24}\\
                                          & \vdots \\
                                          &+\sum_{x \mid D(d)}2\frac{3^{\Omega(p^{2n_2}x)}}{8} + \frac{2^{\Omega(p^{2n_2}x)}}{16}+\frac{1}{24}\\
                                          &\leq 3^{2n_2} (1+2n_2) \cdot \N^{(3)}_{\mathrm{ref}}(L).
 \end{align*}
 and
 \begin{align*}
  \N^{(2)}_{\mathrm{ref}}(K) &= \sum_{x \mid d} \frac{1}{4} \left( 2\frac{2^{\Omega(x)}}{4} + \frac{5}{24} \right)\cdot \left( 2\frac{2^{\Omega(dp^{2n_2}/x)}}{4} + \frac{5}{24} \right)\\ 
                                         &= \sum_{x \mid d} \frac{1}{4} \left( 2\frac{2^{\Omega(px)}}{4} + \frac{5}{24} \right)\cdot \left( 2\frac{2^{\Omega(dp^{2n_2-1}/x)}}{4} + \frac{5}{24} \right) \\
                                         & \vdots \\
                                         &= \sum_{x \mid d} \frac{1}{4} \left( 2\frac{2^{\Omega(p^{2n_2}x)}}{4} + \frac{5}{24} \right)\cdot \left( 2\frac{2^{\Omega(d/x)}}{4} + \frac{5}{24} \right)\\
                                         &\leq 2^{2n_2}(1+2n_2) \cdot \N^{(2)}_{\mathrm{ref}}(L).
 \end{align*}
 Finally note that $\max\{4^{2n_2}, 3^{2n_2} (1+2n_2), 2^{2n_2}(1+2n_2) \} \leq 5103$ for $n_2 \in \{1,2,3\}$.

 Part (a) of this lemma follows from a similar calculation.
 \end{proof}
\end{lem}

\begin{cor}
 Let $L$ be an integral lattice, $p$ an odd prime with $p \nmid \det L$ and $K \in \E_p(L)^{-1}$.
 \begin{itemize}
  \item[(a)] If $K$ is totally-reflective and $\dim L =3$ then 
             \begin{align*}
              81 \cdot \frac{\Nref(L)}{\m(L)} \cdot \left(\left(\frac{1}{1+p^{-1}}\right)^2 p^2 \left(1-p^{-2}\right)\right)^{-1} \geq 1.
             \end{align*}
  \item[(b)] If $K$ is totally-reflective and $\dim L =4$ then 
             \begin{align*}
              5103 \cdot \frac{\Nref(L)}{\m(L)} \cdot \left( \frac{\zeta(4)}{2\zeta(2)^2} \left(\frac{1}{1+p^{-1}}\right)^2 p^3 \left(1-p^{-2}\right)\right)\geq 1.
             \end{align*}
 \end{itemize}
 \begin{proof}
  Combine Lemma $4.6$ and Lemma $4.7$.
 \end{proof}
\end{cor}

\begin{rem}
Since the $p$-term in the above inequalities depends monotonically decreasingly on $p$ (for $p \nmid 2\det L$) and $\Nref(L)/{\m(L)}$ does not depend on $p$ at all, 
it is straightforward to decide when the statement of Corollary $4.8$ is satisfied.
\end{rem}

Given the set of all square free, totally-reflective, primitive genera, one can produce all totally-reflective, primitive genera by using Corollary $4.8$ and Lemma $4.5$.
First, Corollary $4.8$ tells us which (finitely many) primes one has to concider when calculating pre-images under $\E_p$. Then, during the process of repeatedly generating 
lattices $K \in \E_p(L)^{-1}$, Lemma $4.5$ tells us that we can stop and proceed with the next lattice when a non totally-reflective lattice occurs. Eventually this process will terminate
since the number of totally-reflective genera is finite. Furthermore, Lemma $4.5$ implies that every totally-reflective genus will be produced this way.

\section{Reflective Lorentzian lattices and hyperbolic reflection groups}
We want to outline some applications to the theory of reflective Lorentzian lattices and hyperbolic reflection groups. A more detailed investigation will be presented in a (near) future work.

A Lorentzian lattice $E$, that is an integral $\ZZ$-lattice of signature $(n,1)$, is called \emph{reflective} if the Weyl group of $E$ has a finite index in $O(E)$. In terms of reflection groups this 
means that $W(E)$ induces an arithmetic reflection group on the hyperbolic space of dimension $n$ that has a fundamental polyhedron of finite volume. The connection to totally-reflective, 
postive definite lattices is described in the next lemmata.

\begin{lem}
Let $E$ be a square free, reflective lattice of signatrue $(n,1)$ with $n \geq 4$. Then $E$ can be written as $E = {}^{\alpha}\!{\mathbb{H}} \perp L$
with $L$ beeing a square free, totally-reflective, positive definite lattice of dimension $n-1$ and $\alpha \in \ZZ$. If $E$ is strongly square free one can choose $ \alpha = 1$ in the cases $n = 6, n \geq 8$
and $\alpha \in \{1,2\}$ in the cases $n=4,5,7$.
 \begin{proof}
   For $n = 6$ and $n \geq 8$ this follows immediately from \cite{vin81}, Proposition $20$. In cases where this proposition is not applicable one uses elementary manipulations of the genus symbol of $E$
   to see that a (scaled) hyperbolic plane splits $E$. This is possible since the class number of square free Lorentzian lattices is $1$. 
   Vinberg's Lemma implies that the positive definite part is totally-reflective, cf. \cite{vin72}.
 \end{proof}
\end{lem}
The group theoretical version of this lemma reads as follows:

\begin{lem}
Let $W$ be a maximal, arithmetic reflection group on the hyperbolic space of dimension $n \geq 4$. Then $W$ can be written as
 \begin{align*}
  W = W^{+}\left(\mathbb{H} \perp L\right) \quad \text{or} \quad W = W^{+}\left({}^{2}\!{\mathbb{H}} \perp L\right)
 \end{align*}
with $L$ beeing a square free, totally-reflective, positive definite lattice of dimension $n-1$. The $2$ scaled hyperbolic plane is only necessary in the cases $n \in \{4,5\}$.
\begin{proof}
Since $W$ is maximal and arithmetic one can find a Lorentzian lattice $E$ with $W = W(E)$. After repeated use of $\D_p$ and $\E_p$ we can assume that $E$ is strongly square free. The assertion then 
follows from the previous lemma.
\end{proof}
\end{lem}
The following list is a complete classification of all strongly square free, reflective lattices of signature $(5,1)$. As mentioned before, a more detailed proof will be given in a future work.
\begin{thm}
The Lorentzian lattices of signature $(5,1)$ in the table below are reflective. Every square free, reflective Lorentzian lattice of signature $(5,1)$ is isometric to one in the table below.
 \begin{proof}
  Our classification of totally-refelctive lattices provides a list of possible candidates to which we apply Vinberg's algorithm, cf. \cite{vin272}, \S 3.
  In the cases where the algorithm does not terminate, we prove non-reflectivity with the help of a method introduced by Bugaenko in \cite{bug92} and by embedding non-reflective 
  lattices of smaller dimension.
 \end{proof}
\end{thm}

The notation ${}^{\alpha}\!{\mathbb{H}} \perp \text{Genus}$ means that $L$ can be choosen arbitrary within the given $4$-dimensional genus.
The combinatorial structure of the fundamental polyhedron is given as follows:
\begin{itemize}
\item[$r$] = Number of fundamental roots = Number of $4$-dimensional faces,
\item[$f_3$] = Number of $3$-dimensional faces,
\item[$f_2$] = Number of $2$-dimensional faces,
\item[$e$] = Number of edges,
\item[$v$] = Number of vertices,
\item[$c$] = Number of cusps (vertices at infinity).
\end{itemize}

\setlength{\tabcolsep}{2.5mm} 
\begin{center}
 \begin{longtable}{c|c|c|c|c|c|c|c|c}
  No. & \text{$-\det$} & Lattice & $r$ & $f_3$ & $f_2$ & $e$ & $v$ & $c$ \\
  \hline
   & & & & & & & & \\
   1 & 1  & $\mathbb{H} \perp \mathrm{I}(1_{4}^{+4})$  & 6 & 15 & 20 & 15 & 5 & 1 \\
   2 & 2  & $\mathbb{H} \perp \mathrm{I}(2_{5}^{-1})$  & 7 & 20 & 30 & 24 & 8 & 1 \\
   3 & 3  & $\mathbb{H} \perp \mathrm{I}(3^{-1})$  & 7 & 21 & 33 & 27 & 9 & 1 \\
   4 & 3  & $\mathbb{H} \perp \mathrm{I}(3^{+1})$  & 8 & 25 & 40 & 34 & 12 & 1 \\
   5 & 4  & $\mathbb{H} \perp \mathrm{II}(2_{\mathrm{II}}^{-2})$  & 6 & 15 & 20 & 15 & 5 & 1 \\
   6 &4  & $\mathbb{H} \perp \mathrm{I}(2_{2}^{+2})$  & 7 & 21 & 33 & 27 & 8 & 2 \\
   7 & 5  & $\mathbb{H} \perp \mathrm{II}(5^{+1})$  & 6 & 15 & 20 & 15 & 5 & 1 \\
   8 & 5  & $\mathbb{H} \perp \mathrm{I}(5^{+1})$  & 8 & 25 & 40 & 34 & 11 & 2 \\
   9 & 5  & $\mathbb{H} \perp \mathrm{I}(5^{-1})$  & 9 & 32 & 57 & 51 & 18 & 1 \\
  10 &6  & $\mathbb{H} \perp \mathrm{I}(2_{7}^{+1}3^{+1})$  & 9 & 32 & 57 & 51 & 18 & 1 \\
  11 & 6  & $\mathbb{H} \perp \mathrm{I}(2_{5}^{-1}3^{-1})$  & 9 & 31 & 53 & 45 & 14 & 2 \\
  12 & 7  & $\mathbb{H} \perp \mathrm{I}(5^{+1})$  & 11 & 42 & 77 & 70 & 24 & 2 \\
  13 & 8  & $\mathbb{H} \perp \mathrm{I}(2_{3}^{+3})$  & 8 & 25 & 40 & 33 & 10 & 2 \\
  14 & 9  & $\mathbb{H} \perp \mathrm{II}(3^{+2})$   & 7 & 21 & 33 & 27 & 9 & 1 \\
  15 &  9  & $\mathbb{H} \perp \mathrm{I}(3^{+2})$   & 8 & 28 & 50 & 44 & 14 & 2 \\
  16 &  9  & $\mathbb{H} \perp \mathrm{I}(3^{-2})$   & 9 & 32 & 57 & 51 & 16 & 3 \\
  17 & 10 & $\mathbb{H} \perp \mathrm{I}(2_{3}^{-1}5^{-1})$   & 12 & 49 & 94 & 86 & 28 & 3 \\
  18 &  12 & $\mathbb{H} \perp \mathrm{II}(2_{2}^{+2}3^{-1})$   & 7 & 21 & 33 & 27 & 9 & 1 \\
  19 &  12 & $\mathbb{H} \perp \mathrm{II}(2_{6}^{+2}3^{+1})$   & 8 & 25 & 40 & 34 & 12 & 1 \\
  20 & $2^2 \cdot 12$ & ${}^{2}\!{\mathbb{H}} \perp \mathrm{I}(2_{\mathrm{II}}^{+2}3^{-1})$   & 8 & 25 & 40 & 34 & 12 & 1 \\
  21 & $2^2 \cdot 12$ & ${}^{2}\!{\mathbb{H}} \perp \mathrm{I}(2_{\mathrm{II}}^{-2}3^{+1})$   & 7 & 21 & 33 & 27 & 9 & 1 \\
  22 & 12 & $\mathbb{H} \perp \mathrm{I}(2_{0}^{+2}3^{-1})$   & 9 & 33 & 59 & 50 & 14 & 3 \\
  23 &  12 & $\mathbb{H} \perp \mathrm{I}(2_{0}^{+2}3^{+1})$   & 11 & 42 & 77 & 68 & 20 & 4 \\
  24 & 14 & $\mathbb{H} \perp \mathrm{I}(2_{5}^{-1}7^{+1})$   & 15 & 67 & 135 & 127 & 42 & 4 \\
  25 & 15 & $\mathbb{H} \perp \mathrm{I}(3^{+1}5^{+1})$   & 16 & 74 & 153 & 148 & 52 & 3 \\
  26 & 15 & $\mathbb{H} \perp \mathrm{I}(3^{-1}5^{-1})$   & 15 & 66 & 131 & 122 & 40 & 4 \\
  27 & 15 & $\mathbb{H} \perp \mathrm{I}(3^{+1}5^{-1})$   & 12 & 54 & 114 & 113 & 42 & 1 \\
  28 & 18 & $\mathbb{H} \perp \mathrm{I}(2_{3}^{-1}3^{-2})$   & 14 & 62 & 125 & 116 & 36 & 5 \\
  29 & 18 & $\mathbb{H} \perp \mathrm{I}(2_{1}^{+1}3^{+2})$   & 12 & 51 & 101 & 93 & 30 & 3 \\
  30 &  20 & $\mathbb{H} \perp \mathrm{II}(2_{\mathrm{II}}^{+2}5^{+1})$   & 8 & 25 & 40 & 34 & 11 & 2 \\
  31 &  20 & $\mathbb{H} \perp \mathrm{II}(2_{\mathrm{II}}^{-2}5^{-1})$   & 8 & 25 & 40 & 34 & 12 & 1 \\
  32 & 20 & $\mathbb{H} \perp \mathrm{I}(2_{2}^{+2}5^{-1})$   & 16 & 74 & 151 & 140 & 44 & 5 \\
  33 &  21 & $\mathbb{H} \perp \mathrm{II}(3^{-1}7^{+1})$   & 9 & 32 & 57 & 51 & 18 & 1 \\
  34 &  24 & $\mathbb{H} \perp \mathrm{I}(2_{1}^{-3}3^{+1})$   & 11 & 43 & 80 & 71 & 22 & 3 \\
  35 &  24 & $\mathbb{H} \perp \mathrm{I}(2_{5}^{-3}3^{-1})$   & 11 & 43 & 80 & 71 & 22 & 3 \\
  36 & 25 & $\mathbb{H} \perp \mathrm{II}(5^{-2})$   & 9 & 33 & 61 & 57 & 21 & 1 \\
  37 &  25 & $\mathbb{H} \perp \mathrm{I}(5^{-2})$   & 21 & 120 & 282 & 288 & 102 & 5 \\
  38 &   25 & $\mathbb{H} \perp \mathrm{I}(5^{+2})$   & 14 & 67 & 144 & 142 & 46 & 7 \\
  39 &  27 & $\mathbb{H} \perp \mathrm{I}(3^{+3})$   & 9 & 34 & 64 & 58 & 18 & 3 \\
  40 &  27 & $\mathbb{H} \perp \mathrm{I}(3^{-3})$   & 9 & 34 & 64 & 58 & 18 & 3 \\
  41 & 28 & $\mathbb{H} \perp \mathrm{II}(2_{2}^{+2}7^{+1})$   & 11 & 42 & 77 & 70 & 24 & 2 \\
  42 & $2^2 \cdot 28$ & ${}^{2}\!{\mathbb{H}} \perp \mathrm{I}(2_{\mathrm{II}}^{+2}7^{+1})$   & 11 & 42 & 77 & 70 & 24 & 2 \\
  43 &  36 & $\mathbb{H} \perp \mathrm{II}(2_{\mathrm{II}}^{-2}3^{-2})$   & 7 & 21 & 33 & 27 & 8 & 2 \\
  44 &  36 & $\mathbb{H} \perp \mathrm{II}(2_{\mathrm{II}}^{+2}3^{+2})$   & 8 & 28 & 50 & 44 & 14 & 2 \\
  45 &  36 & $\mathbb{H} \perp \mathrm{II}(2_{0}^{+2}3^{+2})$   & 10 & 38 & 69 & 59 & 17 & 3 \\
  46 & $2^2 \cdot 36$ & ${}^{2}\!{\mathbb{H}} \perp \mathrm{I}(2_{\mathrm{II}}^{+2}3^{+2})$   & 10 & 38 & 69 & 59 & 17 & 3 \\
  47 &  36 & $\mathbb{H} \perp \mathrm{I}(2_{2}^{+2}3^{-2})$   & 16 & 83 & 184 & 177 & 54 & 8 \\
  48 &  36 & $\mathbb{H} \perp \mathrm{I}(2_{2}^{+2}3^{+2})$   & 11 & 47 & 94 & 87 & 27 & 4 \\
  49 &  45 & $\mathbb{H} \perp \mathrm{II}(3^{-2}5^{+1})$   & 12 & 50 & 98 & 92 & 30 & 4 \\
  50 &  49 & $\mathbb{H} \perp \mathrm{II}(7^{+2})$   & 16 & 80 & 176 & 176 & 64 & 2 \\
  51 &  54 & $\mathbb{H} \perp \mathrm{I}(2_{3}^{-1}3^{+3})$   & 16 & 75 & 156 & 145 & 42 & 8 \\
  52 &   54 & $\mathbb{H} \perp \mathrm{I}(2_{3}^{-1}3^{-3})$   & 18 & 99 & 230 & 231 & 78 & 6 \\
  53 &  60 & $\mathbb{H} \perp \mathrm{II}(2_{2}^{+2}3^{+1}5^{+1})$   & 16 & 74 & 153 & 148 & 52 & 3 \\
  54 &  60 & $\mathbb{H} \perp \mathrm{II}(2_{2}^{+2}3^{-1}5^{-1})$   & 15 & 66 & 131 & 122 & 40 & 4 \\
  55 &  60 & $\mathbb{H} \perp \mathrm{II}(2_{6}^{+2}3^{+1}5^{-1})$   & 12 & 54 & 114 & 113 & 42 & 1 \\
  56 &  $2^2 \cdot 60$ & ${}^{2}\!{\mathbb{H}} \perp \mathrm{I}(2_{\mathrm{II}}^{+2}3^{+1}5^{+1})$   & 15 & 66 & 131 & 122 & 40 & 4 \\
  57 &  $2^2 \cdot 60$ & ${}^{2}\!{\mathbb{H}} \perp \mathrm{I}(2_{\mathrm{II}}^{+2}3^{-1}5^{-1})$   & 16 & 74 & 153 & 148 & 52 & 3 \\
  58 &  $2^2 \cdot 60$ & ${}^{2}\!{\mathbb{H}} \perp \mathrm{I}(2_{\mathrm{II}}^{+2}3^{-1}5^{+1})$   & 12 & 54 & 114 & 113 & 42 & 1 \\
  59 &  72 & $\mathbb{H} \perp \mathrm{I}(2_{5}^{-3}3^{-2})$   & 24 & 128 & 284 & 274 & 84 & 12 \\
  60 &   72 & $\mathbb{H} \perp \mathrm{I}(2_{7}^{+3}3^{+2})$   & 14 & 66 & 140 & 134 & 44 & 4 \\
  61 &  75 & $\mathbb{H} \perp \mathrm{I}(3^{+1}5^{-2})$   & 86 & 672 & 1788 & 1902 & 660 & 42 \\
  62 &  84 & $\mathbb{H} \perp \mathrm{II}(2_{\mathrm{II}}^{-2}3^{+1}7^{+1})$   & 16 & 74 & 153 & 148 & 52 & 3 \\
  63 &  100 & $\mathbb{H} \perp \mathrm{II}(2_{\mathrm{II}}^{-2}5^{+2})$   & 9 & 33 & 61 & 57 & 19 & 3 \\
  64 &  100 & $\mathbb{H} \perp \mathrm{II}(2_{\mathrm{II}}^{+2}5^{-2})$   & 21 & 120 & 282 & 288 & 102 & 6 \\
  65 & 108 & $\mathbb{H} \perp \mathrm{II}(2_{6}^{+2}3^{-3})$   & 9 & 34 & 64 & 58 & 18 & 3 \\
  66 &108 & $\mathbb{H} \perp \mathrm{II}(2_{6}^{-2}3^{+3})$   & 9 & 34 & 64 & 58 & 18 & 3 \\
  67 & $2^2 \cdot 108$ & ${}^{2}\!{\mathbb{H}} \perp \mathrm{I}(2_{\mathrm{II}}^{+2}3^{-3})$   & 9 & 34 & 64 & 58 & 18 & 3 \\
  68 & $2^2 \cdot 108$ & ${}^{2}\!{\mathbb{H}} \perp \mathrm{I}(2_{\mathrm{II}}^{-2}3^{+3})$   & 9 & 34 & 64 & 58 & 18 & 3 \\
  69 &  108 & $\mathbb{H} \perp \mathrm{I}(2_{6}^{+2}3^{+3})$   & 16 & 85 & 193 & 188 & 56 & 10 \\
  70 &  108 & $\mathbb{H} \perp \mathrm{I}(2_{4}^{-2}3^{-3})$   & 16 & 85 & 193 & 188 & 56 & 10 \\
  71 & 125 & $\mathbb{H} \perp \mathrm{II}(5^{+3})$   & 10 & 40 & 80 & 80 & 30 & 2 \\
  72 & 125 & $\mathbb{H} \perp \mathrm{I}(5^{+3})$   & 20 & 115 & 280 & 295 & 100 & 12 \\
  73 & 180 & $\mathbb{H} \perp \mathrm{II}(2_{\mathrm{II}}^{-2}3^{+2}5^{+1})$   & 16 & 80 & 177 & 178 & 64 & 3 \\
  74 & 196 & $\mathbb{H} \perp \mathrm{II}(2_{\mathrm{II}}^{-2}7^{-2})$   & 15 & 72 & 156 & 159 & 54 & 8 \\
  75 & 216 & $\mathbb{H} \perp \mathrm{I}(2_{7}^{+3}3^{+3})$   & 24 & 138 & 324 & 324 & 104 & 12 \\
  76 & 216 & $\mathbb{H} \perp \mathrm{I}(2_{1}^{-3}3^{-3})$   & 24 & 138 & 324 & 324 & 104 & 12 \\
  77 & 300 & $\mathbb{H} \perp \mathrm{II}(2_{6}^{-2}3^{+1}5^{-2})$   & 86 & 672 & 1788 & 1902 & 660 & 42 \\
  78 & $2^2 \cdot 300$ & ${}^{2}\!{\mathbb{H}} \perp \mathrm{I}(2_{\mathrm{II}}^{-2}3^{-1}5^{-2})$   & 86 & 672 & 1788 & 1902 & 660 & 42 \\
  79 & 500 & $\mathbb{H} \perp \mathrm{II}(2_{\mathrm{II}}^{+2}5^{+3})$   & 20 & 115 & 280 & 295 & 100 & 12 \\
  80 &500 & $\mathbb{H} \perp \mathrm{II}(2_{\mathrm{II}}^{-2}5^{-3})$   & 12 & 56 & 124 & 126 & 44 & 4 \\
 
 \caption{Reflective lattices of siganture $(5,1)$}
 \end{longtable}
\end{center}

\begin{cor}
 Every maximal, arithmetic reflection group on the hyperbolic $5$-space is of the form $W^{+}(E)$ where $E$ is a lattice from the table above.
\end{cor}

\section{Totally-reflective genera}
All discussed calculations were performed using $\mathrm{MAGMA}$, cf. \cite{mag}. 
We used some self-implemented programs to produce all genus symbols of strongly square free genera up to a given bound and determine whether a genus is totally-reflective as well as
an algorithm for calculating pre-images under the Watson-transformation.
Furthermore, we used a program to generate a representative $\ZZ$-lattice for a given genus symbol implemented by Kirschmer \& Lorch, cf. \cite{kl13}.
The (self-implemented) algorithms are available on request.
\subsection{Dimension 3} $1234$ genera, of which $289$ are square free and 52 strongly square free.
\setlength{\tabcolsep}{0.5mm} 
\begin{center}
\begin{longtable}{ccc|ccc|ccc}
  $\det$ & Genus & h & $\det$ & Genus & h & $\det$ & Genus & h \\
  \hline
  & & & & & & & & \\
1&$\mathrm{I}(1_{3}^{+3})$&1&2&$\mathrm{I}(2_{1}^{+1})$&1&3&$\mathrm{I}(3^{-1})$&1\\3&$\mathrm{I}(3^{+1})$&1&4&$\mathrm{I}(2_{2}^{+2})$&1&4&$\mathrm{II}(4_{3}^{-1})$&1\\4&$\mathrm{I}(4_{1}^{+1})$&1&5&$\mathrm{I}(5^{-1})$&1&5&$\mathrm{I}(5^{+1})$&1\\6&$\mathrm{II}(2_{1}^{+1}3^{-1})$&1&6&$\mathrm{I}(2_{1}^{+1}3^{+1})$&1&6&$\mathrm{I}(2_{3}^{-1}3^{-1})$&1\\7&$\mathrm{I}(7^{+1})$&2&8&$\mathrm{I}(2_{1}^{+1}4_{1}^{+1})$&1&9&$\mathrm{I}(3^{-2})$&1\\9&$\mathrm{I}(9^{-1})$&1&9&$\mathrm{I}(9^{+1})$&1&9&$\mathrm{I}(3^{+2})$&1\\
\end{longtable}
\end{center}

\subsection{Dimension 4} $930$ genera, of which $230$ are square free and 88 strongly square free.
\setlength{\tabcolsep}{0.81mm} 
\begin{center}
 \begin{longtable}{ccc|ccc|ccc}
  $\det$ & Genus & h & $\det$ & Genus & h & $\det$ & Genus & h \\
  \hline
  & & & & & & & & \\
1&$\mathrm{I}(1_{4}^{+4})$&1&2&$\mathrm{I}(2_{1}^{+1})$&1&3&$\mathrm{I}(3^{-1})$&1\\3&$\mathrm{I}(3^{+1})$&1&4&$\mathrm{II}(2_{\mathrm{II}}^{-2})$&1&4&$\mathrm{I}(2_{2}^{+2})$&1\\4&$\mathrm{I}(4_{3}^{-1})$&1&4&$\mathrm{I}(4_{1}^{+1})$&1&5&$\mathrm{II}(5^{+1})$&1\\5&$\mathrm{I}(5^{+1})$&1&5&$\mathrm{I}(5^{-1})$&1&6&$\mathrm{I}(2_{1}^{+1}3^{-1})$&2\\6&$\mathrm{I}(2_{1}^{+1}3^{+1})$&1&7&$\mathrm{I}(7^{+1})$&2&8&$\mathrm{I}(2_{3}^{+3})$&1\\8&$\mathrm{II}(2_{1}^{+1}4_{3}^{-1})$&1&8&$\mathrm{I}(2_{1}^{+1}4_{1}^{+1})$&1&9&$\mathrm{II}(3^{+2})$&1\\
\end{longtable}
\end{center}
\bibliographystyle{amsalpha}
\bibliography{References}
 
\end{document}